\documentclass[11pt]{amsart}

\usepackage{lmodern}
\usepackage[colorlinks=true]{hyperref}
\usepackage{amsmath,amssymb,mathrsfs,amsthm, mathtools}
\usepackage[mathscr]{euscript}
\usepackage[utf8]{inputenc}
\usepackage[T1]{fontenc}
\usepackage{mathabx}

\newtheorem{theorem}{Theorem}

\newtheorem{remark}{Remark}

\newcommand{\pv}{\text{p.v.}\int_{\mathbb{S}}}

\def\bbR{{\mathbb R}}
\def\bbS{{\mathbb S^1}}

\def\px{\partial_x}
\def\pt{\partial_t}
\def\unora{\|u(t)\|_{L^2}}
\def\unorb{\|u(t)\|_{\dot{H}^2}}

\title[Nonlocal PDE describing roots of polynomials]{On a nonlocal differential equation describing roots of polynomials under differentiation}
\author[R. Granero-Belinch\'{o}n]{Rafael Granero-Belinch\'{o}n}
\email{rafael.granero@unican.es}
\address{Departamento  de  Matem\'aticas,  Estad\'istica  y  Computaci\'on,  Universidad  de Cantabria.  Avda.  Los  Castros  s/n,  Santander,  Spain.}

\date{\today}

\allowdisplaybreaks
\begin{document}

\begin{abstract}
In this work we study the nonlocal transport equation derived recently by Steinerberger when studying how the distribution of roots of a polynomial behaves under iterated differentation of the function. In particular, we study the well-posedness of the equation, establish some qualitative properties of the solution and give conditions ensuring the global existence of both weak and strong solutions. Finally, we present a link between the equation obtained by Steinerberger and a one-dimensional model of the surface quasi-geostrophic equation used by Chae, C\'ordoba, C\'ordoba \& Fontelos.
\end{abstract}

\maketitle

\section{Introduction and main results}
In this paper, we consider the following one-dimensional nonlinear transport equation
\begin{equation}\label{eq:1}
\pt u+\px\arctan\left(\frac{Hu}{u}\right)=0\qquad (x,t)\text{ on }\bbS\times[0,T],
\end{equation}
where 
$$
H u(x)=\frac{1}{2\pi}\pv \frac{u(y)}{\tan\left(\frac{x-y}{2}\right)}dy,
$$
is the Hilbert transform. The previous equation needs to be supplemented with the initial data
\begin{equation}\label{eq:2}
u(x,0)=u_0(x).
\end{equation}

This equation has been derived by S. Steinerberger \cite{St2018} when studying how the distribution of roots behaves under iterated differentation. In particular, the purpose of this work is to study the properties of the transport equation \eqref{eq:1}. The study of such nonlocal and nonlinear one-dimensional equations is a wide research area with a large literature. For other similar equations and related results we refer to \cite{MR2456270, MR2680191, MR812343, MR2179734, MR3171771, MR3621816, MR3492724, MR2464570, MR3542965, MR2397459, MR2439488,LaLe2015}.

In this paper we prove the following results:

\begin{theorem}\label{teo1}
Let $0< u_0\in H^2(\bbS)$ be the initial data. Then there exists a time $0<T\leq\infty$, $T=T(\|u_0\|_{H^2},\min_x u_0(x))$ and a unique positive solution to \eqref{eq:1}
$$
0<u\in C([0,T],H^2(\bbS)).
$$
Furthermore, this solution verifies the following properties: 
\begin{itemize}
\item
\begin{align}\label{identity:1}
\unora^2+\int_0^t\mathscr{D}(s) ds= \|u_0\|_{L^2}^2.
\end{align}
where
$$
\mathscr{D}=\frac{1}{16\pi}\int_\bbR\pv\frac{u(x)-u(y)}{\sin\left(\frac{x-y}{2}\right)^2}\log\left(\frac{u(x)^2+(Hu(x))^2}{u(y)^2+(Hu(y))^2}\right)dxdy,
$$
\item  $\|u(t)\|_{L^1}=\|u_0\|_{L^1},$
\item if $u_0(x)$ is even, $u(x,t)$ remains even,
\item Maximum principle: $\max_x u(x,t)\leq \max_x u_0(x),$
\item Minimum principle: $\min_x u_0(x)\leq \min_x u(x,t).$
\end{itemize}
\end{theorem}
\begin{remark}
We remark that, for an arbitrary $u(x,t)$, we are not able to give a sign to $\mathscr{D}$ (compare with \cite{MR3281135,granero2017fractional}). In other words, we are not able to show whether the $L^2$ norm decays.
\end{remark}
\begin{remark}
The following scaling is invariant for the equation $u_\lambda(x,t)=\lambda^{\alpha}u(\lambda x,\lambda^{1-\alpha} t)\;\forall\,\alpha\in\mathbb{R}.$
\end{remark}

Under certain restrictions we can ensure that the solution is global:
\begin{theorem}\label{teo2}
Let $0< u_0\in H^2(\bbS)$ be the initial data and denote 
$$
\langle u_0\rangle=\frac{1}{2\pi}\int_\bbS u_0(x)dx.
$$
There exists $0<\mathcal{C}$ such that if
$$
\frac{\left\|u_0-\langle u_0\rangle\right\|_{A^1}}{\langle u_0\rangle}\leq \mathcal{C},
$$
then the solution (from Theorem \ref{teo1}) is global and satisfies
$$
\|u(t)-\langle u_0\rangle\|_{A^1}\leq \|u_0-\langle u_0\rangle\|_{A^1}e^{-\delta t}
$$
for certain $0<\delta(\langle u_0\rangle)$ small enough.
\end{theorem}

\begin{remark}
The explicit lower bound $0.17<\mathcal{C}$ is obtained as a byproduct.
\end{remark}

Finally, we study the existence of weak solutions \emph{i.e.} solutions that satisfy the equation in the following sense:
\begin{equation*}
-\int_{0}^T\int_{\mathbb{S}^1} u(x,s)\pt \phi(x,s)+\arctan\left(\frac{Hu(x,s)}{u(x,s)}\right)\px \phi dxds=\int_{\mathbb{S}^1} u_0(x) \phi(x,0) dx,
\end{equation*}
for all test functions $\phi(x,t)\in C^\infty(\mathbb{S}^1\times[0,T))$.

In that regards, we prove the global existence of weak solution for initial data satisfying certain size conditions in a scale invariant space (with respect to the scaling of the equation):
\begin{theorem}\label{teo3}
Let $0< u_0\in A^0(\bbS)$ be the initial data and denote 
$$
\langle u_0\rangle=\frac{1}{2\pi}\int_\bbS u_0(x)dx.
$$
There exists $0<\tilde{\mathcal{C}}$ such that if
$$
\frac{\left\|u_0-\langle u_0\rangle\right\|_{A^0}}{\langle u_0\rangle}\leq \tilde{\mathcal{C}},
$$
then there exists at least one global weak solution
$$
u\in L^\infty([0,T]\times \mathbb{S}^1)\cap L^2(0,T;H^{0.5}),\forall\;0<T<\infty
$$ 
and this solution satisfies
$$
\|u(t)-\langle u_0\rangle\|_{L^\infty}\leq \|u_0-\langle u_0\rangle\|_{A^0}e^{-\delta t}
$$
for certain $0<\delta(\langle u_0\rangle)$ small enough. 
\end{theorem}

\begin{remark}
The explicit lower bound $0.24<\tilde{\mathcal{C}}$ is obtained as a byproduct.
\end{remark}

\begin{remark}
Similar results can be proved for the porous medium equation (see also \cite{stan2016finite, stan2014finite})
$$
\pt u+\px \left(\frac{Hu}{u^m}\right)=0,\;m\in\mathbb{N}.
$$
\end{remark}

The rest of the paper is devoted to the proofs of each results (sections \ref{sec:1}-\ref{sec:3}) and the link between \eqref{eq:1} and the equation
\begin{equation}\label{eq:hmodel}
\pt g+\Lambda g=\px\left(gHg\right),
\end{equation}
(see section \ref{sec:4}). We would like to remark that \eqref{eq:hmodel} was proposed as a one-dimensional model of the 2D surface quasi-geostrophic equation by Chae, C\'ordoba, C\'ordoba \& Fontelos \cite{MR2141858} (see also the papers by Matsuno \cite{matsuno1991linearization} and Baker, Li \& Morlet \cite{baker1996analytic}). 

\subsection*{Notation}
We denote
$$
\Lambda u=H\px u(x)=\frac{1}{4\pi}\pv\frac{u(x)-u(x-y)}{\sin^2(y/2)}dy.
$$
We define the $L^2$ based Sobolev spaces 
$$
H^s=\left\{u(x)=\sum_{n\in\mathbb{Z}} \hat{u}(n)e^{inx}\text{ with }\sum_{n\in\mathbb{Z}} |n|^{2s}|\hat{u}(n)|^2<\infty \right\},
$$
with norm $\|u\|_{H^s}=\|\Lambda^s u\|_{L^2}.$
Similarly, we recall the definition of the Wiener spaces
$$
A^s=\left\{u(x)=\sum_{n\in\mathbb{Z}} \hat{u}(n)e^{inx}\text{ with }\sum_{n\in\mathbb{Z}} |n|^{s}|\hat{u}(n)|<\infty \right\}.
$$
with norm $\|u\|_{A^s}=\|\widehat{\Lambda^s u}\|_{\ell^1}.$
\section{Proof of Theorem \ref{teo1}}\label{sec:1}
\subsection*{Well-posedness} 
The existence will follow using the energy method \cite{majda2002vorticity} once the appropriate \emph{a priori} estimates are obtained. We define the energy
\begin{equation}\label{ineq:0}
\mathscr{E}(t)=\frac{1}{\min_x u(x,t)}+\|u(t)\|_{H^2}.
\end{equation}
We have to proof an inequality of the type
$$
\frac{d}{dt}\mathscr{E}(t)\leq C(1+\mathscr{E}(t))^p,
$$
for certain $C$ and $p$.

To estimate the first term in the energy we use a pointwise argument (see \cite{MR2084005,MR2989434,MR3238306, MR3469062} for more details). The solution has at least a minimum:
$$
m(t)=\min_ x u(x,t)=u(\underline{x}_t,t).
$$
Because of the positivity of the initial data, we have that $m(0)>0$. Following the argument in \cite{MR2084005,MR3238306, MR3469062}, we have that
$$
\frac{d}{dt}m(t)=\pt u(\underline{x}_t,t)=-\frac{m(t) \Lambda u(\underline{x}_t)}{m(t)^2 + (Hu(\underline{x}_t))^2}\;\text{ a.e.}.
$$
Then,
\begin{equation}\label{ineq:1}
\frac{d}{dt}\frac{1}{\min_x u(x,t)}=-\frac{\pt u(\underline{x}_t,t)}{m(t)^2}\leq C\frac{\|u\|_{H^2}}{m(t)^3}\leq C(\mathscr{E}(t))^4
\end{equation}

For the sake of brevity we only provide with the estimates for the higher order terms (being the rest of the terms straightforward). We take 2 derivatives of the equation and test against $\px^2u$. We find that
\begin{align*}
\frac{d}{dt}\unorb^2&=-\frac{(\px^2u \Lambda u+u \Lambda \px^2 u )\px^2u}{u^2 + (Hu)^2}+\frac{u \Lambda u (2u\px^2u+2Hu\px\Lambda u)\px^2u}{(u^2 + (Hu)^2)^2}\\
&\quad+\frac{(\px\Lambda u \px u+Hu \px^3 u)\px^2 u}{u^2 + (Hu)^2}-\frac{Hu \px u(2u\px^2u+2Hu\Lambda\px u)\px^2 u}{(u^2 + (Hu)^2)^2}+\text{l.o.t.}\\
&=I_1+I_2+I_3+I_4+\text{l.o.t.}.
\end{align*}
Using that
$$
\left\|\frac{1}{u^2+(Hu)^2}\right\|_{L^\infty}\leq \left\|\frac{1}{u^2}\right\|_{L^\infty}\leq \frac{1}{m(t)^2}\leq \mathscr{E}(t)^2,
$$
we have that
\begin{equation}\label{ineq:2}
\text{l.o.t.}\leq C(\mathscr{E}(t))^8\unorb.
\end{equation}
We recall the C\'ordoba-C\'ordoba inequality \cite{MR2032097} 
$$
\theta\Lambda \theta\geq \frac{1}{2}\Lambda(\theta^2),
$$
to find that
$$
\int_\bbS \frac{u \Lambda \px^2 u \px^2u}{u^2 + (Hu)^2} dx\geq\frac{1}{2}\int_\bbS H\px\left(\frac{u}{u^2 + (Hu)^2}\right)(\px^2 u)^2dx
$$
\begin{align*}
I_1&\leq \frac{\unorb^2\|u\|_{A^1}}{m(t)^2}-\frac{1}{2}\int_\bbS H\px\left(\frac{u}{u^2 + (Hu)^2}\right)(\px^2 u)^2dx\\
&\leq \frac{\unorb^2\|u\|_{A^1}}{m(t)^2}+\frac{1}{2}\frac{\unorb^2\|u\|_{A^1}}{m(t)^2}+\frac{\unorb^2\|u\|_{A^0}^2\|u\|_{A^1}}{m(t)^4}\\
&\leq C(1+\mathscr{E}(t))^8\unorb.
\end{align*}
A similar use of H\"{o}lder inequality and Sobolev embedding leads to
$$
I_2+I_4\leq C(1+\mathscr{E}(t))^8\unorb.
$$
Finally, we observe that an integration by parts allow us to conclude
$$
I_3\leq C(1+\mathscr{E}(t))^8\unorb.
$$
Then, we obtain that
\begin{equation}\label{ineq:3}
\frac{1}{2}\frac{d}{dt}\unorb^2\leq C(1+\mathscr{E}(t))^8\unorb.
\end{equation}
Thus, collecting \eqref{ineq:1},\eqref{ineq:2} and \eqref{ineq:3}, we conclude te desired inequality \eqref{ineq:0}

\subsection*{An identity for the evolution of the $L^2$ norm} 
We test \eqref{eq:1} against $u$. We find that
\begin{align*}
\frac{d}{dt}\unora^2&=\int_\bbS-\frac{u \Lambda u u}{u^2 + (Hu)^2}+\frac{Hu \px u u}{u^2 + (Hu)^2}dx.
\end{align*}
We compute
\begin{align*}
\mathscr{D}&=-\frac{1}{2}\int_\bbS Hu\px\log\left(u^2+(Hu)^2\right)dx\\
&=-\int_\bbS Hu\frac{u\px u+Hu\Lambda u}{u^2+(Hu)^2}dx\\
&=-\int_\bbS \frac{u\px u Hu}{u^2+(Hu)^2}-\frac{u^2\Lambda u}{u^2+(Hu)^2}+\Lambda u dx\\
&=-\int_\bbS \frac{u\px u Hu}{u^2+(Hu)^2}-\frac{u^2\Lambda u}{u^2+(Hu)^2}dx.
\end{align*}
As a consequence,
\begin{align*}
\frac{d}{dt}\unora^2&=-\mathscr{D}.
\end{align*}
Furthermore,
\begin{align*}
\mathscr{D}&=\frac{1}{2}\int_\bbS \Lambda u\log\left(u^2+(Hu)^2\right)dx\\
&=\frac{1}{8\pi}\int_\bbS\pv\frac{u(x)-u(x-y)}{\sin(y/2)^2}\log\left(u(x)^2+(Hu(x))^2\right)dxdy\\
&=\frac{1}{8\pi}\int_\bbS\pv\frac{u(x)-u(y)}{\sin((x-y)/2)^2}\log\left(u(x)^2+(Hu(x))^2\right)dxdy\\
&=\frac{1}{8\pi}\int_\bbS\pv\frac{u(y)-u(x)}{\sin((x-y)/2)^2}\log\left(u(y)^2+(Hu(y))^2\right)dxdy.
\end{align*}
Then, we have identity \eqref{identity:1}. This concludes with the existence part. The uniqueness follow from a standard contradiction argument.

\subsection*{Propagation of the $L^1$ norm}
Once the solution remains positive, the $L^1$ norm is preserved due to the divergence form of the equation.

\subsection*{Propagation of the even symmetry} This is a straightforward consequence of the fact that the Hilbert transform $H$ maps even functions into odd functions.

\subsection*{Maximum principle}
We define
$$
M(t)=\max_ x u(x,t)=u(\overline{x}_t,t)
$$
Then (see \cite{granero2014global,MR3238306, MR3469062} for more details) we have that
$$
\frac{d}{dt}M(t)=\pt u(\overline{x}_t,t)\;\text{ a.e.},
$$
Then
$$
\frac{d}{dt}M(t)+\frac{M(t) \Lambda u(\overline{x}_t)}{M(t)^2 + (Hu(\overline{x}_t))^2}=0,
$$
We observe that 
$$
\Lambda u(\overline{x}_t)\geq0.
$$
Thus, using $0\leq M(t)$, we obtain that
$$
M(t)\leq M(0).
$$
\subsection*{Minimum principle}
With the previous definition for $m(t)$, we have that
$$
\frac{d}{dt}m(t)+\frac{m(t) \Lambda u(\underline{x}_t)}{m(t)^2 + (Hu(\underline{x}_t))^2}=0.
$$
Thus, using $\Lambda u(\underline{x}_t)\leq 0$, we find that
$$
0\leq m(0)\leq m(t).
$$

\section{Proof of Theorem \ref{teo2}}
The proof of this Theorem follows the approach in \cite{granero2018global}. We define the new variable
$$
v(x,t)=u(x,t)-\langle u_0\rangle.
$$
This variable quantifies the difference between the steady state $u_\infty=\langle u_0\rangle$ and $u$. The idea of the theorem is first to linearize around the stady state $\langle u_0\rangle$. Secondly, we obtain an inequality of the form 
$$
\frac{d}{dt}\|v(t)\|_{A^1}+\frac{\|v(t)\|_{A^2}}{\langle u_0\rangle}\leq \mathcal{F}\left(\frac{\|v(t)\|_{A^1}}{\langle u_0\rangle}\right)\|v(t)\|_{A^2},
$$
with $\mathcal{F}(0)=0$ and $\mathcal{F}$ continuous. We observe that this inequality guarantees $v(t)\rightarrow0$ in $A^1$ for small enough $\|v_0\|_{A^1}/\langle u_0\rangle$. 

In what follows we assume that
$$
r=\frac{\|v\|_{A^1}}{\langle u_0\rangle}<\frac{1}{2},
$$
so that
$$
\frac{\|v\|_{A^1}}{\langle u_0\rangle-\|v\|_{A^1}}=\frac{r}{1-r}<1.
$$
Since we have the following Poincar\'e type inequality
$$
\|v\|_{A^s}\leq \|v\|_{A^{r}},\; \forall 0\leq s<r,
$$
we observe that
$$
\bigg{|}\frac{Hu}{u}\bigg{|}=\bigg{|}\frac{Hv}{u}\bigg{|}\leq \frac{\|v\|_{A^0}}{\langle u_0\rangle-\|v\|_{A^0}}\leq \frac{\|v\|_{A^1}}{\langle u_0\rangle-\|v\|_{A^1}}<1.
$$
As a consequence, we can expand the nonlinearity as a power series 
$$
\arctan\left(\frac{Hu}{u}\right)=\sum_{n\in \mathbb{Z}^+\cup\{0\}}\frac{(-1)^n}{1+2n}\left(\frac{Hu}{u}\right)^{1+2n},
$$
so
$$
\pt u=-\sum_{n\in \mathbb{Z}^+\cup\{0\}}(-1)^n\left(\frac{Hu}{u}\right)^{2n}\left(\frac{\Lambda u}{u}-\frac{Hu \px u}{u^2}\right).
$$
In the new variable, this latter equation reads
\begin{align*}
\pt v&=-\sum_{n\in \mathbb{Z}^+}(-1)^n\left(\frac{Hv}{v+\langle u_0\rangle}\right)^{2n}\left(\frac{\Lambda v}{v+\langle u_0\rangle}-\frac{Hv \px v}{(v+\langle u_0\rangle)^2}\right)\\
&\quad-\left(\frac{\Lambda v}{v+\langle u_0\rangle}-\frac{Hv \px v}{(v+\langle u_0\rangle)^2}\right).
\end{align*}
We recall the following Taylor series
$$
\frac{1}{\langle u_0\rangle +v}=\frac{1}{\langle u_0\rangle}+\frac{1}{\langle u_0\rangle}\sum_{n\in \mathbb{Z}^+}(-1)^{n}\left(\frac{v}{\langle u_0\rangle}\right)^{n},
$$
$$
\frac{1}{(\langle u_0\rangle +v)^2}=\frac{1}{\langle u_0\rangle^2}+\frac{1}{\langle u_0\rangle^2}\sum_{n\in \mathbb{Z}^+}(-1)^{n}(1+n)\left(\frac{v}{\langle u_0\rangle}\right)^{n}.
$$
We define
\begin{align*}
\mathfrak{S}_1&=\frac{Hv}{\langle u_0\rangle}+\frac{Hv}{\langle u_0\rangle}\sum_{m\in \mathbb{Z}^+}(-1)^{m}\left(\frac{v}{\langle u_0\rangle}\right)^{m}\\
\mathfrak{S}_2&=\frac{Hv \px v}{\langle u_0\rangle^2}+\frac{Hv \px v}{\langle u_0\rangle^2}\sum_{m\in \mathbb{Z}^+}(-1)^{m}(1+m)\left(\frac{v}{\langle u_0\rangle}\right)^{m}\\
\mathfrak{S}_3&=\frac{\Lambda v}{\langle u_0\rangle}\sum_{m\in \mathbb{Z}^+}(-1)^{m}\left(\frac{v}{\langle u_0\rangle}\right)^{m}.
\end{align*}
Using the previous Taylor series together with the previous definitions, we find that
\begin{align}
\pt v+\frac{\Lambda v}{\langle u_0\rangle}&=-\sum_{n\in \mathbb{Z}^+}(-1)^n\left(\mathfrak{S}_1\right)^{2n}\left(\frac{\Lambda v}{\langle u_0\rangle}+\mathfrak{S}_3\right)\nonumber\\
&\quad+\sum_{n\in \mathbb{Z}^+}(-1)^n\left(\mathfrak{S}_1\right)^{2n}\mathfrak{S}_2-\mathfrak{S}_3+\mathfrak{S}_2\label{eq:v}.
\end{align}
We take a derivative of \eqref{eq:v} to obtain that 
\begin{align}
\pt \px v+\frac{\Lambda \px v}{\langle u_0\rangle}&=\sum_{n\in \mathbb{Z}^+}(-1)^{n+1}2n\left(\mathfrak{S}_1\right)^{2n-1}\px \mathfrak{S}_1\left(\frac{\Lambda v}{\langle u_0\rangle}+\mathfrak{S}_3\right)\nonumber\\
&\quad-\sum_{n\in \mathbb{Z}^+}(-1)^n\left(\mathfrak{S}_1\right)^{2n}\left(\frac{\Lambda \px v}{\langle u_0\rangle}+\px \mathfrak{S}_3\right)\nonumber\\
&\quad+\sum_{n\in \mathbb{Z}^+}(-1)^n2n\left(\mathfrak{S}_1\right)^{2n-1}\px\mathfrak{S}_1\mathfrak{S}_2+\sum_{n\in \mathbb{Z}^+}(-1)^n\left(\mathfrak{S}_1\right)^{2n}\px \mathfrak{S}_2\nonumber\\
&\quad-\px\mathfrak{S}_3+\px\mathfrak{S}_2\label{eq:v2}
\end{align}
We want to estimate 
$$
\|v\|_{A^1}=\|\px v\|_{A^0}.
$$
To do that we first observe that $A^0$ is an algebra, thus,
\begin{align*}
\|\mathfrak{S}_1^{2n}\|_{A^0}\leq\|\mathfrak{S}_1\|_{A^0}^{2n}\leq \left(\frac{\|v\|_{A^0}}{\langle u_0\rangle}+\frac{\|v\|_{A^0}}{\langle u_0\rangle}\sum_{m\in \mathbb{Z}^+}\left(\frac{\|v\|_{A^0}}{\langle u_0\rangle}\right)^{m}\right)^{2n}.
\end{align*}
Summing up the series, we find the estimate
\begin{align*}
\|\mathfrak{S}_1^{2n}\|_{A^0}\leq \left(\frac{\|v\|_{A^0}}{\langle u_0\rangle-\|v\|_{A^0}}\right)^{2n}.
\end{align*}
Similarly,
\begin{align*}
\|\mathfrak{S}_1^{2n-1}\|_{A^0}&\leq \left(\frac{\|v\|_{A^0}}{\langle u_0\rangle-\|v\|_{A^0}}\right)^{2n-1}\\
\|\mathfrak{S}_2\|_{A^0}&\leq \frac{\|v\|_{A^0} \|v\|_{A^1}}{\langle u_0\rangle^2}+\frac{\|v\|_{A^0} \|v\|_{A^1}}{\langle u_0\rangle^2}\sum_{m\in \mathbb{Z}^+}(1+m)\left(\frac{\|v\|_{A^0}}{\langle u_0\rangle}\right)^{m}\leq \frac{\|v\|_{A^0} \|v\|_{A^1}}{(\langle u_0\rangle -\|v\|_{A^0})^2}\\
\|\mathfrak{S}_3\|_{A^0}&\leq\frac{\|v\|_{A^1}}{\langle u_0\rangle}\sum_{m\in \mathbb{Z}^+}\left(\frac{\|v\|_{A^0}}{\langle u_0\rangle}\right)^{m}\leq \frac{\|v\|_{A^1}\|v\|_{A^0}}{\langle u_0\rangle (\langle u_0\rangle -\|v\|_{A^0})}\\
\|\px \mathfrak{S}_1\|_{A^0}&\leq \frac{\|v\|_{A^1}}{\langle u_0\rangle-\|v\|_{A^0}}+\frac{\|v\|_{A^1}\|v\|_{A^0}}{(\langle u_0\rangle-\|v\|_{A^0})^2}\\
\|\px \mathfrak{S}_2\|_{A^0}&\leq \frac{\|v\|_{A^1}^2+\|v\|_{A^0} \|v\|_{A^2}}{(\langle u_0\rangle -\|v\|_{A^0})^2}+\frac{2\|v\|_{A^0}\|v\|_{A^1}^2}{(\langle u_0\rangle -\|v\|_{A^0})^3}\\
\|\px \mathfrak{S}_3\|_{A^0}&\leq \frac{\|v\|_{A^2}\|v\|_{A^0}}{\langle u_0\rangle (\langle u_0\rangle -\|v\|_{A^0})}+\frac{\|v\|_{A^1}^2}{(\langle u_0\rangle-\|v\|_{A^0})^2}.
\end{align*}
We obtain that
\begin{align}
\frac{d}{dt}\|v\|_{A^1}+\frac{\|v\|_{A^2}}{\langle u_0\rangle}&\leq \frac{\|v\|_{A^2}}{\langle u_0\rangle}\bigg{\{}\sum_{n\in \mathbb{Z}^+}2n\left(\frac{r}{1-r}\right)^{2n-1}\left[\frac{r}{1-r}+\frac{r^2}{(1-r)^2}\right]\frac{1}{1-r}\nonumber\\
&\quad+\sum_{n\in \mathbb{Z}^+}\left(\frac{r}{1-r}\right)^{2n}\left(1+\frac{r}{1-r}+\frac{r}{(1-r)^2}\right)\nonumber\\
&\quad+\sum_{n\in \mathbb{Z}^+}2n\left(\frac{r}{1-r}\right)^{2n-1}\left[\frac{r}{1-r}+\frac{r^2}{(1-r)^2}\right]\frac{r}{(1-r)^2}\nonumber\\
&\quad+2\sum_{n\in \mathbb{Z}^+}\left(\frac{r}{1-r}\right)^{2n}\left(\frac{r}{(1-r)^2}+\frac{r^2}{(1-r)^3}\right)\nonumber\\
&\quad+\frac{r}{1-r}+\frac{r}{(1-r)^2}\bigg{\}}.\label{eq:v3}
\end{align}
Using
$$
\frac{z^2}{1-z^2}=\sum_{n\in\mathbb{Z}^+}z^{2n},
$$
we find
\begin{align}
\frac{d}{dt}\|v\|_{A^1}+\frac{\|v\|_{A^2}}{\langle u_0\rangle}&\leq \frac{\|v\|_{A^2}}{\langle u_0\rangle}\bigg{\{}\frac{2\frac{r}{1-r}}{\left(1-\left(\frac{r}{1-r}\right)^2\right)^2}\left[\frac{r}{1-r}+\frac{r^2}{(1-r)^2}\right]\frac{1}{1-r}\nonumber\\
&\quad+\frac{\left(\frac{r}{1-r}\right)^2}{1-\left(\frac{r}{1-r}\right)^2}\left(1+\frac{r}{1-r}+\frac{r}{(1-r)^2}\right)\nonumber\\
&\quad+\frac{2\frac{r}{1-r}}{\left(1-\left(\frac{r}{1-r}\right)^2\right)^2}\left[\frac{r}{1-r}+\frac{r^2}{(1-r)^2}\right]\frac{r}{(1-r)^2}\nonumber\\
&\quad+2\frac{\left(\frac{r}{1-r}\right)^2}{1-\left(\frac{r}{1-r}\right)^2}\left(\frac{r}{(1-r)^2}+\frac{r^2}{(1-r)^3}\right)\nonumber\\
&\quad+\frac{r}{1-r}+\frac{r}{(1-r)^2}\bigg{\}}.\label{eq:v4}
\end{align}
Finally, we can simplify the previous expression and find that
\begin{align*}
\mathcal{F}(r)&=\frac{2\frac{r}{1-r}}{\left(1-\left(\frac{r}{1-r}\right)^2\right)^2}\left(\frac{r}{(1-r)^3}+\frac{r}{(1-r)^4}\right)+\frac{\left(\frac{r}{1-r}\right)^2}{1-2r}\nonumber\\
&\quad+2\frac{\left(\frac{r}{1-r}\right)^2}{1-\left(\frac{r}{1-r}\right)^2}\frac{1}{(1-r)^3}+\frac{r}{1-r}+\frac{r}{(1-r)^2}.
\end{align*}
We observe that $\mathcal{F}$ is a continuous function in a neighborhood of $r=0$ and satisfies $\mathcal{F}(0)=0$. Thus, it exists $0<\mathcal{C}$ such that $\mathcal{F}(\mathcal{C})<1$. We finally observe that if $\|v_0\|_{A^1}/\langle u_0\rangle < \mathcal{C}$ this condition propagates in time and ensures the following bound
$$
\|v(t)\|_{A^1}\leq \|v_0\|_{A^1}e^{-\delta t},
$$
for small enough $0<\delta\ll1$. This last inequality together with a close inspection of the energy estimates in Theorem 1 lead to the following inequality
$$
\frac{d}{dt}\|u\|_{H^2}^2\leq C(u_0)\|u\|_{H^2}^2,
$$
and then we conclude the global bound for the $H^2$ norm using Gronwall's inequality.

\section{Proof of Theorem \ref{teo3}}\label{sec:3}
In this section we prove the existence of global weak solutions for certain initial data satisfying appropriate size restriction in the space $A^0$. We emphasize that this space is scale invariant with respect to the scaling of the equation. First we obtain a priori estimates, then we consider a vanishing viscosity approximation and prove the convergence of the approximate solutions.

\subsection*{\emph{A priori} estimates} 
Following the previous ideas, the first nonlinear term contributes with
\begin{align*}
\|NL_1\|_{A^0}&=\left\|\sum_{n\in \mathbb{Z}^+}(-1)^n\left(\mathfrak{S}_1\right)^{2n}\left(\frac{\Lambda v}{\langle u_0\rangle}+\frac{\Lambda v}{\langle u_0\rangle}\sum_{m\in \mathbb{Z}^+}(-1)^{m}\left(\frac{v}{\langle u_0\rangle}\right)^{m}\right)\right\|_{A^0}\\
&\leq \sum_{n\in \mathbb{Z}^+}\left(\frac{\|v\|_{A^0}}{\langle u_0\rangle-\|v\|_{A^0}}\right)^{2n}\left(\frac{\|v\|_{A^1}}{\langle u_0\rangle-\|v\|_{A^0}}\right)\\
&\leq\left(\frac{\|v\|_{A^1}}{\langle u_0\rangle-\|v\|_{A^0}}\right) \left(\frac{1}{1-\left(\frac{\|v\|_{A^0}}{\langle u_0\rangle-\|v\|_{A^0}}\right)^2}-1\right)\\
&\leq\left(\frac{\|v\|_{A^1}}{\langle u_0\rangle-\|v\|_{A^0}}\right) \left(\frac{\left(\frac{\|v\|_{A^0}}{\langle u_0\rangle-\|v\|_{A^0}}\right)^2}{1-\left(\frac{\|v\|_{A^0}}{\langle u_0\rangle-\|v\|_{A^0}}\right)^2}\right).
\end{align*}
The second nonlinear term can be estimated as
\begin{align*}
\|NL_2\|_{A^0}&=\left\|\sum_{n\in \mathbb{Z}^+}(-1)^n\left(\mathfrak{S}_1\right)^{2n}\left(\frac{Hv \px v}{\langle u_0\rangle^2}+\frac{Hv \px v}{\langle u_0\rangle^2}\sum_{m\in \mathbb{Z}^+}(-1)^{m}(1+m)\left(\frac{v}{\langle u_0\rangle}\right)^{m}\right)\right\|_{A^0}\\
&\leq \sum_{n\in \mathbb{Z}^+}\left(\frac{\|v\|_{A^0}}{\langle u_0\rangle-\|v\|_{A^0}}\right)^{2n}\left(\frac{\|v\|_{A^0}\|v\|_{A^1}}{(\langle u_0\rangle-\|v\|_{A^0})^2}\right)\\
&\leq\left(\frac{\|v\|_{A^0}\|v\|_{A^1}}{(\langle u_0\rangle-\|v\|_{A^0})^2}\right) \left(\frac{\left(\frac{\|v\|_{A^0}}{\langle u_0\rangle-\|v\|_{A^0}}\right)^2}{1-\left(\frac{\|v\|_{A^0}}{\langle u_0\rangle-\|v\|_{A^0}}\right)^2}\right).
\end{align*}
Finally, we find that
\begin{align*}
\|NL_3\|_{A^0}&=\left\|\Lambda v\frac{1}{\langle u_0\rangle}\sum_{n\in \mathbb{Z}^+}(-1)^{n}\left(\frac{v}{\langle u_0\rangle}\right)^{n}\right\|_{A^0}\\
&\leq \frac{\|v\|_{A^1}}{\langle u_0\rangle}\sum_{n\in \mathbb{Z}^+}\left(\frac{\|v\|_{A^0}}{\langle u_0\rangle}\right)^{n}\\
&\leq \|v\|_{A^1}\left(\frac{1}{\langle u_0\rangle-\|v\|_{A^0}}-\frac{1}{\langle u_0\rangle}\right),
\end{align*}
\begin{align*}
\|NL_4\|_{A^0}&=\left\|\left(\frac{Hv \px v}{\langle u_0\rangle^2}+\frac{Hv \px v}{\langle u_0\rangle^2}\sum_{n\in \mathbb{Z}^+}(-1)^{n}(1+n)\left(\frac{v}{\langle u_0\rangle}\right)^{n}\right)\right\|_{A^0}\\
&\leq \frac{\|v\|_{A^0}\|v\|_{A^1}}{\langle u_0\rangle^2}+\frac{\|v\|_{A^0}\|v\|_{A^1}}{\langle u_0\rangle^2}\sum_{n\in \mathbb{Z}^+}(1+n)\left(\frac{\|v\|_{A^0}}{\langle u_0\rangle}\right)^{n}\\
&\leq \frac{\|v\|_{A^0}\|v\|_{A^1}}{(\langle u_0\rangle-\|v\|_{A^0})^2}.
\end{align*}
We define
$$
s=\frac{\|v\|_{A^0}}{\langle u_0\rangle}.
$$
Collecting the previous estimates, we find that
\begin{align*}
\frac{d}{dt}\|v\|_{A^0}+\frac{\|v\|_{A^1}}{\langle u_0\rangle}&\leq \frac{\|v\|_{A^1}}{\langle u_0\rangle}\bigg{[}\frac{s}{(1-s)^2}+\frac{s}{1-s}+\left(\frac{s}{(1-s)^2}\right) \left(\frac{\left(\frac{s}{1-s}\right)^2}{1-\left(\frac{s}{1-s}\right)^2}\right)\\
&\quad+\left(\frac{1}{1-s}\right) \left(\frac{\left(\frac{s}{1-s}\right)^2}{1-\left(\frac{s}{1-s}\right)^2}\right)\bigg{]}
\end{align*}
Using the hypotheses on $\tilde{\mathcal{C}}$, we conclude that
$$
s\leq \tilde{\mathcal{C}}
$$
implies 
$$
\frac{d}{dt}\|v\|_{A^0}+\delta\|v\|_{A^1}\leq 0,
$$
and that, thank to a Poincar\'e-type inequality, leads to 
$$
\|v(t)\|_{L^\infty}\leq \|v_0\|_{A^0}e^{-\delta t}.
$$
Furthermore, the solution also enjoys the following parabolic gain of regularity
$$
\int_0^t\|v(s)\|_{H^{0.5}}^2ds\leq \int_0^t\|v(s)\|_{A^1}ds\sup_{s}\|v(s)\|_{L^1}\leq \int_0^t\|v(s)\|_{A^1}ds\|v_0\|_{A^0}2\pi.
$$
\subsection*{Approximated solutions} 
To construct the approximate solutions, we consider the following vanishing viscosity approximated problem \begin{equation}\label{eq:1eps}
\pt u^ \varepsilon+\px\arctan\left(\frac{Hu^ \varepsilon}{u^ \varepsilon}\right)=\varepsilon\px^2 u^ \varepsilon\qquad (x,t)\text{ on }\bbS\times[0,T],
\end{equation}
with a mollified initial data
$$
u^\varepsilon(x,0)=\mathscr{M}_\varepsilon*u_0(x).
$$
The corresponding approximate solution exists globally and remains smooth.

\subsection*{Compactness}
Fix $0<T<\infty$. We have that $u^\varepsilon$ is uniformly bounded in 
$$
L^\infty(0,T;A^0)\cap L^2(0,T;H^{0.5}).
$$
This implies weak-* convergence
$$
u^\varepsilon\overset{\ast}{\rightharpoonup} u,
$$
in
$$
L^\infty([0,T]\times\mathbb{S}^1),
$$
and weak convergence
$$
u^\varepsilon\rightharpoonup u,
$$
in
$$
L^2(0,T;H^{0.5}(\mathbb{S}^1)).
$$
Furthermore, $\pt u^\epsilon$ is uniformly bounded in
$$
L^2(0,T;H^{-1.5}).
$$
A standard application of Aubin-Lions Theorem \cite{temam2001navier} ensures the strong convergence (after maybe taking a subsequence)
$$
u^\varepsilon\rightarrow u,\;Hu^\varepsilon\rightarrow Hu
$$
in 
$$
L^2(0,T;L^2).
$$
Taking another subsequence if necessary, we obtain that 
$$
u^\varepsilon(x,t)\rightarrow u(x,t)\;\text{ a.e in }\mathbb{S}^1\times[0,T]
$$
In particular, we conclude the lower bound
$$
\min_{x} u_0(x)\leq u(x,t)\text{ a.e in }\mathbb{S}^1\times[0,T].
$$

\subsection*{Passing to the limit}
Being the other terms linear, we only have to take into consideration the convergence of 
$$
J=\int_0^T\int_{\mathbb{S}^1}\left(\arctan\left(\frac{Hu^ \varepsilon}{u^ \varepsilon}\right)-\arctan\left(\frac{Hu}{u}\right)\right)\px\phi dxds.
$$
We have that
$$
J\leq \int_0^T\int_{\mathbb{S}^1}\left|\frac{Hu^ \varepsilon}{u^ \varepsilon}-\frac{Hu}{u}\right||\px\phi |dxds.
$$
Using the lower bounds for $u$ and $u^\epsilon$ together with H\"older inequality, we conclude that
$$
J\rightarrow 0.
$$
This concludes the proof of the existence of a global weak solution $u$.

\section{Link between (\ref{eq:1}) and (\ref{eq:hmodel})}\label{sec:4}
We now look for a solution of \eqref{eq:1} having the following form
$$
u(x,t)=\langle u_0\rangle + \varepsilon\sum_{j=0}^\infty \varepsilon^j f^{(j)}(x,t),
$$
(here $\varepsilon$ can be though as the displacement from the homogeneous state $\langle u_0\rangle$). The idea is to truncate the series up to certain order, two say,
$$
f(x,t)=\varepsilon f^{(0)}(x,t)+\varepsilon^2 f^{(1)}(x,t),
$$
and see what $f$ solves. In this way we will obtain that (up to $O(\varepsilon^3)$), $f$ solves \eqref{eq:hmodel}. A similar approach has been used in the study of free boundary problems for incompressible fluids (see \cite{cheng2018rigorous, granero2018asymptotic, granero2018asymptotic2} and the references therein). First, we observe that \eqref{eq:1} can be equivalently written as
\begin{equation}\label{eq:3}
\pt u+\frac{u \Lambda u - Hu \px u}{u^2 + (Hu)^2}=0\qquad (x,t)\text{ on }\bbS\times[0,T].
\end{equation}
Thus,
$$
\pt u\left(\langle u_0\rangle^2+2(u-\langle u_0\rangle)\langle u_0\rangle+(u-\langle u_0\rangle)^2 + (Hu)^2\right)+u \Lambda u - Hu \px u=0.
$$
Forcing the previous ansatz and matching the powers of $\varepsilon$, we find that $f^{(0)}$ solves
$$
\pt f^{(0)}+\frac{\Lambda  f^{(0)}}{\langle u_0\rangle}=0.
$$
Similarly, $f^{(1)}$ solves
$$
\pt f^{(1)}\langle u_0\rangle^2+2\pt f^{(0)}f^{(0)}\langle u_0\rangle+\langle u_0\rangle \Lambda f^{(1)}+f^{(0)}\Lambda f^{(0)} - Hf^{(0)} \px f^{(0)}=0.
$$
Thus, substituting $\langle u_0\rangle\pt f^{(0)}$ by $-\Lambda  f^{(0)}$, we find that $f$ solves
$$
\pt f+\frac{1}{\langle u_0\rangle}\Lambda f - \frac{1}{\langle u_0\rangle^2}\px (Hf f)=O(\varepsilon^3).
$$
Thus, neglecting the $O(\varepsilon^3)$ terms we find that 
$$
g(x,t)=\frac{f(x,t\langle u_0\rangle)}{\langle u_0\rangle}
$$ 
solves \eqref{eq:hmodel}.

\subsection*{Acknowledgements}R. G-B has been funded by the grant MTM2017-89976-P from the Spanish government.

\bibliographystyle{plain}

\end{document}